\documentclass[11pt]{article}
\usepackage{}
\usepackage{amsfonts}
\usepackage{amssymb}
\usepackage{amsmath,amssymb,graphicx}
\usepackage{epsfig}
\newtheorem{theorem}{Theorem}[section]
\newtheorem{theorema}{Theorem}

\newtheorem{lemma}[theorem]{Lemma}

\newcommand{\kz}{\kappa(z)}
\newcommand{\nn}{\nonumber}

\newcommand{\be}{\begin{eqnarray}}
\newcommand{\ee}{\end{eqnarray}}
\newcommand{\ba}{\begin{array}}
\newcommand{\ea}{\end{array}}
\newcommand{\ben}{\begin{eqnarray*}}
\newcommand{\een}{\end{eqnarray*}}

\newcommand{\pa}{\partial}

\begin{document}

\hoffset -34pt

\title{{\large \bf Asymptotic properties of the hyperbolic metric on the sphere with three conical singularities}
\thanks{\mbox{Keywords. Conical singularities, hyperbolic metrics.}}}
\author{\normalsize  Tanran Zhang }
\date{}
\maketitle \baselineskip 21pt
\noindent

\begin{minipage}{138mm}
\renewcommand{\baselinestretch}{1} \normalsize
\begin{abstract}
{The explicit formula for the hyperbolic metric $\lambda_{\alpha,\,\beta,\,\gamma}(z)|dz|$ on the thrice-punctured sphere $\mathbb{P} \backslash \{z_1,\,z_2,\,z_3\}$ with singularities of order $\alpha,\,\beta,\,\gamma \leq 1$ with $\alpha+\beta+\gamma>2$ at $z_1,\,z_2,\,z_3$ was given by Kraus, Roth and Sugawa in \cite{Rothhyper}. In this paper we investigate the asymptotic properties of the higher order derivatives of $\lambda_{\alpha,\,\beta,\,\gamma}(z)$ near the singularity and give some more precise description for the asymptotic behavior.  }
\end{abstract}
\end{minipage}\\
\\
\\
\renewcommand{\baselinestretch}{1} \normalsize
\section{Introduction}

From the viewpoint of the theory of partial differential equations and functions, the hyperbolic metric, also called the $\mathrm{Poincar\acute{e}}$ metric, plays an important role in metric spaces. It facilitates describing the hyperbolic geometry on some domains in different ways. Since the Gaussian curvature of the hyperbolic metric is a negative constant, it can be regarded as the extremal metric of a class of regular conformal metrics with strictly negative Gaussian curvature functions. This kind of conformal metric is more general, and it was discussed by Heins in \cite{Heins}, Kraus, Roth and Ruscheweyh in \cite{Roth2007}.

Equipped with the hyperbolic metric, a punctured domain is more complicated than a simply connected domain. An elementary case for the punctured domain is the hyperbolic metric $\lambda_{\mathbb{D}^*}(z)|dz|$ on the (once-)punctured unit disk $\mathbb{D} \backslash \{0\}$, which is defined by
$$\lambda_{\mathbb{D}^*}(z)|dz|=\frac{|dz|}{2|z|\log(1/|z|)}$$
with the Gaussian curvature $-4$. It is induced by the hyperbolic metric
$$\lambda_{\mathbb{D}}(z)|dz|=\frac {|dz|}{1-{|z|}^2}$$
on the unit disk $\mathbb{D}$. The expression for the hyperbolic metric on the twice-punctured unit disk is not achieved yet, and there are only a few estimates for the density function, pre-Schwarzian and Schwarzian, see \cite{Hempel1979, Hempel1989, Smith1989} for details.
\par In the thrice-punctured sphere $\mathbb{P} \backslash \{z_1,\,z_2,\,z_3\}$ with singularities of order $\alpha,\,\beta,\,\gamma \leq 1$ at $z_1,\,z_2,\,z_3$, if $\alpha+\beta+\gamma>2$, the hyperbolic metric $\lambda_{\alpha,\,\beta,\,\gamma}(z)|dz|$ can be expressed in terms of special functions in $\mathbb{C} \backslash \{0,\,1\}$, see \cite{Rothhyper}. Kraus, Roth and Sugawa used the Liouville equation
\be\label{liouville equation}
\Delta u=4e^{2u}
\ee
to obtain the explicit formula of $\lambda_{\alpha,\,\beta,\,\gamma}(z)$. For equation (\ref{liouville equation}), Liouville proved in \cite{Liouville1} that, in any disk $D$ contained in the punctured unit disk $\mathbb{D}\backslash\{0\}$ every solution $u$ to \eqref{liouville equation} can be written as
\be \label{solution}
u(z)=\log\frac{|f'(z)|}{1-|f(z)|^2},
\ee
where $f$ is a holomorphic function in $D$. Kraus, Roth and Sugawa at first obtained the expression of function $f$ in equation (\ref{solution}) using the hypergeometric differential equation
$$z(1-z)w''(z)+[\alpha-(\alpha+\beta-1)z]w'(z)-\frac{(\alpha+\beta-\gamma)(\alpha+\beta+\gamma-2)}{4}w(z)=0,$$
and then gave the explicit formula for $\lambda_{\alpha,\,\beta,\,\gamma}(z)$. Note that only the special equation (\ref{liouville equation}) is involved here since $\lambda_{\alpha,\,\beta,\,\gamma}(z)$ has the Gaussian curvature $-4$. We concern the estimate for the derivatives of $\lambda_{\alpha,\,\beta,\,\gamma}(z)$ near the origin and give a stronger result then the estimates in \cite{Zhang1}.
\par Then we discuss the so-called Minda-type theorems. In 1997, Minda \cite{Mindametric} studied the behavior of the hyperbolic metric in a neighborhood of a puncture on the plane domain using the uniformization theorem for up to second order derivatives. His results can be extended to higher order derivatives of a conformal metric with negative curvatures on an arbitrary hyperbolic region, see \cite{Zhang1}. However, near the origin, if the order $\alpha<1$, this kind of limits may not exist. We prove that for the hyperbolic metric $\lambda_{\alpha,\,\beta,\,\gamma}(z)$, the limits in Minda-type always exist, and give the recurrence formula for them.

\section {Preliminaries }
\setcounter{equation}{0}

For complex numbers $a$, $b$, $c$ with $c\neq 0$, $-1$, $-2$, $\ldots$, the Gaussian hypergeometric function is defined as
$$F(a,b,c;z)=\sum^{\infty}_{n=0}\frac{(a)_n(b)_n}{(c)_n}\frac{z^n}{n!}, \quad \quad |z|<1,$$ where $(a)_n$ is the Pochhammer symbol, namely, $(a)_0=1$ and
\be\label{bracket}
(a)_n=a(a+1)\cdots (a+n-1) \nn
\ee
for $n=1,\,2,\,3, \ldots \ $. It is continued analytically to the slit plane $\mathbb{C} \backslash [1, +\infty)$. Its derivative is given by
\be \label{geometric derivative}
\frac{d}{d z}F(a,b,c;z)=\frac{ab}{c}F(a+1,b+1,c+1;z).
\ee
We can immediately obtain
\be \label{geometric higher derivative}
\frac{d^n}{d z^n}F(a,b,c;z)=\frac{(a)_n(b)_n}{(c)_n}F(a+n,b+n,c+n;z).
\ee
We have
\be \label{asym terms}
F(a,b,c;z)=\frac{\Gamma(c)\Gamma(c-a-b)}{\Gamma(c-a)\Gamma(c-b)}F(a,b,a+b-c+1;1-z) \qquad \qquad \qquad \qquad \quad \nn \\
+(1-z)^{c-a-b} \frac{\Gamma(c)\Gamma(a+b-c)}{\Gamma(a)\Gamma(b)}F(c-a,c-b,c-a-b+1;1-z)
\ee
for $|\arg(1-z)|<\pi$, where $\Gamma(z)$ is the gamma function, see 15.3.6 in \cite{handbook1965}. Each term of \eqref{asym terms} has a pole when $c=a+b \pm n$, $n =0,\,1,\,2,\,\ldots$, and this case is covered by
\be \label{covering case}
&&F(a,b,a+b+n; z) \nn \\
&=&\frac{\Gamma(n)\Gamma(a+b+n)}{\Gamma(a+n)\Gamma(b+n)} \sum^{n-1}_{j=0}\frac{(a)_j(b)_j}{j!(1-n)_j}(1-z)^j \nn \\
&&-\frac{\Gamma(a+b+n)}{\Gamma(a)\Gamma(b)}(z-1)^n  \sum^{\infty}_{j=0}\frac{(a+n)_j(b+n)_j}{j!(j+n)!}(1-z)^j[\log(1-z) \nn \\
&&\qquad -\Psi(j+1)-\Psi(j+n+1)+\Psi(a+j+n)+\Psi(b+j+n)],
\ee
for $|\arg(1-z)|<\pi$, $|1-z|<1$, where $\Psi(z)={\Gamma'(z)}/{\Gamma(z)}$ is the digamma function, see 15.3.11 in \cite{handbook1965}, and we take the convention that $\sum_{j=a}^b=0$ if $b < a$ here and after. The behavior of the hypergeometric function near $z=1$ satisfies
\be \label{behavior of hypergeometric}
\left\{
\begin{array}{l} \vspace*{2mm}
\displaystyle F(a,b,c;1)=\frac{\Gamma(c)\Gamma(c-a-b)}{\Gamma(c-a)\Gamma(c-b)}, \ \textrm{if}\ a+b<c, \\ \vspace*{2mm}
\displaystyle F(a,b,a+b;z)=\frac{1}{B(a,b)} \left( \log {\frac{1}{1-z}}+R(a,b)\right)(1+\textit{O}(1-z)), \\ \vspace*{2mm}
F(a,b,c;z)=(1-z)^{c-a-b}F(c-a,c-b,c;z), \ \textrm{if}\ a+b>c.\vspace*{1mm}\end{array} \right. \ee
Here
\be \label{beta}
B(a,b)=\frac{\Gamma(a)\Gamma(b)}{\Gamma(a+b)}
\ee
is the beta function and
\be \label{rab}
R(a,b)=2\Psi(1)-\Psi( a)-\Psi( b)
\ee
with
$\Psi(x)={\Gamma'(x)}/{\Gamma(x)}$ being the digamma function. The asymptotic formula in (\ref{behavior of hypergeometric}) for the case $a+b=c$ is due to Ramanujan, see \cite{handbook1965, Anderson2009}.\\

In the domain $G \subseteq \mathbb{C}$, every positive, upper semi-continuous function $\lambda: G \rightarrow (0, + \infty)$ induces a conformal metric on $G$. We denote the metric by $\lambda(z)|dz|$, see \cite{Heins,Roth2007}. We $\lambda(z)|dz|$ the metric and $\lambda(z)$ the density to avoid any ambiguity. A conformal metric $\lambda(z)|dz|$ on a domain $G\subseteq\mathbb{C}$ is said to be regular, if its density $\lambda(z)$ is positive and twice continuously differentiable on $G$, i.e. $\lambda(z)>0$ and $\lambda(z) \in C^2(G)$. For a domain $G \subseteq \mathbb{C}$ equipped with a conformal metric $\lambda(z)|dz|$, the distance function associating two points $z, \zeta \in G$ is defined by
\be\label{distance}
\rho_{\lambda}(z,\zeta):=\inf_{\iota} \int_{\iota}\lambda(z)|dz|,
\ee
where the infimum is taken over all rectifiable paths $\iota$ in $G$ joining $z$ and $\zeta$. We call $(G, \rho_{\lambda})$ a metric space. The metric $\lambda(z)|dz|$ is said to be complete on $G$ if $(G, \,\rho_{\lambda})$ is a complete metric space. The Gaussian curvature $\kappa_{\lambda}(z)$ of the regular conformal metric $\lambda(z)|dz|$ is defined by
$$\kappa_{\lambda}(z)=-\frac{\Delta\log\lambda(z)}{\lambda(z)^2},$$
where $\Delta$ denotes the Laplace operator. For the definition of the Gaussian curvature in more general case, see \cite{zhang2}.
\par The basic property of Gaussian curvature is its conformal invariance. That means, given a (regular) conformal metric $\lambda(z)|dz|$ on a domain $G$ and a holomorphic mapping $f: \Omega\rightarrow G$ on a Riemann surface $\Omega$, the pullback
\be
f^{*}\lambda(w)|dw|:=\lambda(f(w))|f'(w)||dw| \nn
\ee
is still a (regular) conformal metric on on $\Omega \backslash\{\mbox{critical points of}\ f\}$ with Gaussian curvature
\be
\kappa_{f^{*}{\lambda}}(w)=\kappa_{\lambda}(f(w)). \nn
\ee
Here $\Omega \backslash\{\mbox{critical points}$ $\mbox{of}\ f\}$ is a punctured domain, the critical points of $f$ are the source of punctures. If the neighborhood of a puncture carries some special structure as given below, we say this puncture is a singularity.
\par Let $\mathbb{P}$ denote the Riemann sphere $\mathbb{C}\cup \{\infty\}$ and let $\Omega\subseteq\mathbb{P}$ be a subdomain. For a point $p \in \Omega$, let $z$ be local coordinates such that $z(p)=0$. We say a conformal metric $\lambda(z)|dz|$ on the punctured domain $\Omega^*:=\Omega \backslash \{p\}$ has a conical singularity of order $\alpha\leq 1$ at the point $p$, if, in local coordinates $z$,
\be \label{singularity}
\log\lambda(z)=\left\{
\begin{array}{ll}
-\alpha\log|z|+v(z) & \mbox{if\ }\ \alpha < 1 \\
-\log|z|-\log\log(1/|z|)+w(z)&\mbox{if\ }\ \alpha=1,\end{array} \right. \ee
where $v(z), w(z) = \textit {O}(1)$ as $z(p)\rightarrow 0$ with $\textit {O}$ and $\textit {o}$ being the Landau symbols throughout our study. For $u(z):=\log \lambda(z)$, the order $\alpha$ of $\lambda(z)|dz|$ at the conical singularity $p$ is again the order of $u(z)$ at the conical singularity $\log p$. We call the point $p$ a corner of order $\alpha$ if $\alpha< 1$ and a cusp if $\alpha=1$. It is evident that the cusp is the limit case of the corner.

The hyperbolic metric $\lambda_{\Omega}(z)dz$ on a domain $\Omega$ is a complete metric with some negative constant Gaussian curvature, here we take the constant to be  $-4$. By \eqref{distance}, the hyperbolic distance between $z,\zeta \in \Omega$ is
$$\texttt{d}_{\Omega}(z,\zeta):=\inf_{\iota} \int_{\iota}\lambda_{\Omega}(z)|dz|$$
and the infimum is always attained; the hyperbolic line passing through $z$ and $\zeta$ is the path for which the infimum is attained. From the conformal invariancy of Gaussian curvature we know that, the hyperbolic metric $\lambda_{\Omega}(z)|dz|$ on any domain $\Omega$ induces a hyperbolic metric on some domain which is conformally equivalent to $\Omega$. The hyperbolic metric is a kind of metric of special interest because it is the unique maximal conformal metric in the sense of conformal invariance, see \cite{Ahlforslemma, Heins}. The following result gives the explicit formula of the hyperbolic metric on twice-punctured plane $\mathbb{C} \backslash \{0,\,1\}$. The terminology \emph{generalized hyperbolic metric} is motivated by the fact that if all singularities are cusps, then we can get back the standard hyperbolic metric on the punctured sphere $\mathbb{P} \backslash \{z_1,\,\ldots,\,z_n\}$, see \cite{Rothhyper}.
\begin{theorema} $(\mathrm{[10]})$ \label{explicit formula}
\textsl{Let $0<\alpha,\ \beta <1$ and $0<\gamma \leq 1$ such that $\alpha+\beta+\gamma>2$. Then the generalized hyperbolic density on the thrice-punctured sphere $\mathbb{P}\backslash \{0,\,1,\,\infty\}$ of orders $\alpha,\, \beta,\,\gamma$ at $0$, $1$, $\infty$, respectively, can be expressed by
\be
&& \lambda_{\alpha,\,\beta,\,\gamma}(z) \nonumber\\
&=& \frac 1{|z|^{\alpha}|1-z|^{\beta}}\cdot\frac{K_3}{K_1|\varphi_1(z)|^2
+K_2|\varphi_2(z)|^2+2Re(\varphi_1(z)\varphi_2(\bar{z}))} \label{lambda(z)1} \\
&=&\frac{1}{|z|^{\alpha}|1-z|^{\beta}} \cdot \frac{\delta (1-\alpha)}{|\varphi_1(z)|^2-{\delta^2}|1-z|^{2-2\alpha}|{\varphi_3}(z)|^2} \label{lambda(z)2}
\ee
in the twice-punctured plane $\mathbb{C}\backslash \{0,\,1\}$, where
\be
K_1:=-\frac{\Gamma( c-a)\Gamma( c- b)}
{\Gamma( c)\Gamma( c-a- b)},\ \ \
K_2:=-\frac{\Gamma(a+1- c)\Gamma( b+1- c)}
{\Gamma(1- c)\Gamma(a+ b+1- c)}, \label{K1K2} \\
K_3:=\sqrt{\frac{\sin(\pi a)\sin(\pi b)}
{\sin(\pi( c- a))\sin(\pi( c- b))}}
\cdot\frac{\Gamma( a+ b+1- c)\Gamma( c)}{\Gamma( a)\Gamma( b)}\qquad  \nn
\ee
and
$$\varphi_1(z)=F( a, b, c;z),\ \ \ \varphi_2(z)=F( a, b, a+ b- c+1;1-z),$$
$$\varphi_3(z)=F( a-c+1, b-c+1, 2-c;z),$$
with
\be \label{abc}
a=\frac{\alpha+\beta-\gamma}{2},\ b=\frac{\alpha+\beta+\gamma-2}{2},\  c=\alpha;
\ee
\be \label{delta}
\delta=\frac{\Gamma( c)}{\Gamma(2- c)}
\left(\frac{\Gamma(1- a)\Gamma(1- b)\Gamma( a+1- c)
\Gamma( b+1- c)}{\Gamma( a)\Gamma( b)\Gamma( c- a)
\Gamma( c- b)} \right)^{1/2}.
\ee
The Gaussian curvature of $\lambda(z)$ defined by (\ref{lambda(z)1}) and (\ref{lambda(z)2}) is $-4$. Note that $\varphi_1$ and $\varphi_3$ are analytic in $\mathbb{C}\backslash [1,\,+\infty)$, $\varphi_2$ is analytic in $\mathbb{C}\backslash (-\infty,\,0]$.}
\end{theorema}

Expressions (\ref{lambda(z)1}) and (\ref{lambda(z)2}) are equal to each other. Denote $\log \lambda(z):=\log \lambda_{\alpha,\,\beta,\,\gamma}(z)$ for short, and
$$\partial^n:=\frac{\partial ^n}{\partial z^n},\ \bar{\partial}^n:=\frac {\partial^n}{\partial\bar{z}^n}$$
for $n\geq1$. The following theorem is a general estimate for $\log \lambda(z)$ near the singularities.
\begin{theorema} $(\mathrm{[15]})$ \label{general pu}
\textsl{For $\lambda(z)$ as in \eqref{lambda(z)1} with order $\alpha\in(0,\,1]$, let $u(z)=\log\lambda(z)$. Then for $m,\,n\geq 1$, \vspace*{2mm} \\
(i) $ \displaystyle \ \lim_{z\rightarrow0}z^n\partial^nu(z)=\frac {\alpha}2(-1)^n(n-1)!=\lim_{z\rightarrow0}\bar{z}^n\bar{\partial}^nu(z),$ \vspace*{2mm} \\
(ii) $ \displaystyle \ \lim_{z\rightarrow 0}\bar{z}^m z^n\bar{\partial}^m\partial^nu(z)=0.$}
\end{theorema}

We can estimate the higher order derivatives of a conformal density function $\lambda(z)$ directly. The following result is of Minda-type.
\begin{theorema} $(\mathrm{[15]})$ \label{coro}
\textsl{ Let $\lambda(z)|dz|$ be a regular conformal metric on a domain $\Omega\subseteq\mathbb{C}$ with an isolated singularity at $z=p$. Suppose that the curvature $\kappa :\Omega\rightarrow \mathbb{R}$ has a H\"{o}lder continuous extension to $\Omega\cup\{p\}$ such that $\kappa(p)< 0$ and the order of $\log\lambda$ is $\alpha=1$ at $z=p$. Then \vspace*{2mm} \\
(i)$\ \displaystyle \lim_{z\rightarrow p} (z-p)|z-p|\log(1/|z-p|)\lambda_{z}(z)=-{\frac 1 {2\sqrt{-\kappa(p)}}}, \hspace*{\fill} $ \vspace*{2mm}\\
(ii)$\ \displaystyle \ \lim_{z\rightarrow p}(z-p)^2|z-p|\log(1/|z-p|)\lambda_{zz}(z)=\displaystyle {\frac 3{4\sqrt{-\kappa(p)}}}, \hspace*{\fill} $ \vspace*{2mm}\\
(iii)$\ \displaystyle \ \lim_{z\rightarrow p}|z-p|^3\log(1/|z-p|)\lambda_{z\bar{z}}(z)=\displaystyle {\frac 1 {4\sqrt{-\kappa(p)}}}. \hspace*{\fill} $ }
\end{theorema}

Theorem \ref{coro} was given only for the order $\alpha=1$. When the order $\alpha<1$, the analogous limit
\be \label{lza}
\lim_{z\rightarrow0}|z|^{\alpha}\lambda(z)
\ee
does not necessarily exist. But for the hyperbolic density $\lambda_{\alpha,\,\beta,\,\gamma}(z)$, if $0<\alpha <1$, expression (\ref{lambda(z)2}) shows that the limit (\ref{lza}) exists. The following theorem is due to Kraus, Roth and Sugawa in \cite{Rothhyper}. They did not give the explicit formula of (\ref{lza}), but it is easy to deduce that from Corollary 4.4 in their paper.
\begin{theorema} \label{limit for alpha less than 1}
\textsl{For the hyperbolic density $\lambda_{\alpha, \beta, \gamma}$ given in (\ref{lambda(z)2}), if $0<\alpha<1$, then we have
\be \label{l0}
\lim_{z \rightarrow0} |z|^{\alpha} \lambda(z)=\frac{\delta}{1-\delta^2}(1-\alpha)
\ee
where $\delta$ is as in (\ref{delta}), $a$, $b$ and $c$ are as in (\ref{abc}).}
\end{theorema}

\section{Case $0<\alpha <1$}
\setcounter{equation}{0}

In this section we consider the hyperbolic metric when the order $0<\alpha <1$. We again let  $\lambda(z):=\lambda_{\alpha, \beta, \gamma}(z)$. For the hyperbolic density function $\lambda(z)$, we can only consider the asymptotic behavior near the origin. By the expression of $\lambda(z)$, we know that the singularity $z=1$ is the same as the origin. As for the infinity, we can change the coordinates by a conformal function, say, $z \mapsto 1/z$, to map $\infty$ onto $0$. But some calculation is involved, so it is convenient to consider the case near the origin. In expression \eqref{lambda(z)1}, for orders $0<\alpha,\,\beta<1$ and $0<\gamma\leq 1$, the real parameters $\alpha, \beta, \gamma$ given by condition \eqref{abc} satisfy
$$-\frac{1}{2}<a<1,\ -1<b<\frac{1}{2},\ 0<c<1.$$
At first we give a lemma for future use.
\begin{lemma} \label{Mz lemma}
\textsl{In the expression for $\lambda(z)$ as in \eqref{lambda(z)1} with order $\alpha\in(0,\,1)$, let
\be \label{Mz}
M(z):&=&K_1|\varphi_1(z)|^2+K_2|\varphi_2(z)|^2+2Re\left(\varphi_1(z)\varphi_2(\bar{z})\right) \nn \\
&=&(K_1\varphi_1(\bar{z})+\varphi_2(\bar{z}))\varphi_1(z)+(K_2\varphi_2(\bar{z})
+\varphi_1(\bar{z}))\varphi_2(z).
\ee
Then for $a,\,b$ and $c$ are defined in \eqref{abc}, $K_1$ and $K_2$ are defined in \eqref{K1K2}, \vspace*{2mm}\\
(1)$ \displaystyle \ \lim_{z\rightarrow 0}\pa M(z)=\frac{ab}{c}\left(K_1-\frac{1}{K_2}\right)\ $ for $\displaystyle \ 0<\alpha<{1}/{2}$, \vspace*{1mm}\\
(2)$ \displaystyle \ \pa M(z)=2ab\left(K_1-\frac{1}{K_2}\right)+2K_2 \left(\frac{\Gamma(c)\Gamma(a+b-c+1)}{\Gamma(a)\Gamma(b)}\right)^2 \frac{\bar{z}}{|z|} +O\left(|z|^{\frac{1}{2}}\right)$  for $\alpha={1}/{2}$,\vspace*{1mm}\\
(3)$ \displaystyle \ \lim_{z\rightarrow 0}z^{n} |z|^{2\alpha-2}\pa^{n} M(z)=\frac{(-1)^{n-1}(c)_{n-1}K_2}{1-c}\left(\frac{\Gamma(c)\Gamma(a+b-c+1)}{\Gamma(a)\Gamma(b)}\right)^2$ \vspace*{1mm}\\
for $n\geq 2$ if $\ 0<\alpha\leq {1}/{2}$ and $n\geq 1$ if $\ {1}/{2}<\alpha<1$, \\
(4)\vspace*{2mm} $ \displaystyle \lim_{z\rightarrow 0}\bar{z}^{m}z^{n} |z|^{2\alpha-2}\bar{\pa}^{m}\pa^{n} M(z)=(-1)^{n+m}(c)_{n-1}(c)_{m-1}K_2 \left(\frac{\Gamma(c)\Gamma(a+b-c+1)}{\Gamma(a)\Gamma(b)}\right)^2$ \\
for $m,\,n \geq 1$. }
\end{lemma}
\textbf{Remark. }Case $(2)$ can be expressed by $\pa M(z)=O(1)$. It is easy to see that there is no non-vanishing limit such as in $(3)$ holds for $n=1$ and $\alpha={1}/{2}$, even if it is multiplied by a power of $z/\bar{z}$.\vspace*{3mm} \\
\textbf{Proof of Theorem \ref{Mz lemma}. }Since $\varphi_1(z),\,\varphi_2(z)$ are analytic in $\mathbb{C} \backslash [1,\,+\infty)$, $\mathbb{C} \backslash (-\infty,\,0]$ respectively, then we have $\overline{\pa^n \varphi_1(z)}=\bar{\pa}^n \left(\varphi_1(\bar{z})\right)$ for $z\in \mathbb{C} \backslash [1,\,+\infty)$, $\overline{\pa^n \varphi_2(z)}=\bar{\pa}^n \left(\varphi_2(\bar{z})\right)$ for $z \in \mathbb{C} \backslash (-\infty,\,0]$.
For limit (1), we have
\ben
\pa M(z)=(K_1\varphi_1(\bar{z})+\varphi_2(\bar{z}))\pa \varphi_1(z)+(K_2\varphi_2(\bar{z})
+\varphi_1(\bar{z}))\pa \varphi_2(z).
\een
From properties \eqref{geometric higher derivative},
\be \label{phi11} \displaystyle
\pa \varphi_1(0)=\frac{ab}{c}, \ee
and from \eqref{behavior of hypergeometric},
\be \label{phi1 and phi2}
\displaystyle \varphi_2(0)=F(a,b,a+b-c+1;1)=-\frac 1 {K_2},\ \varphi_1(0)=1,
\ee
provided that $a+b<a+b-c+1$, so
\be \label{former adder}
K_1\varphi_1(0)+\varphi_2(0)=K_1-{K_2}^{-1}.
\ee
Now we consider the term $(K_2\varphi_2(\bar{z})+\varphi_1(\bar{z}))\pa \varphi_2(z)$, which satisfies $$\lim_{z\rightarrow 0}(K_2\varphi_2(z)+\varphi_1(z))=0.$$
Note that
\ben
&&\varphi_2(z)=F(a,b,a+b-c+1; 1-z) \\
&=&\frac{\Gamma(a+b-c+1)\Gamma(1-c)}{\Gamma(b-c+1)\Gamma(a-c+1)}F(a, b, c; z) \\
&& +z^{1-c}\frac{\Gamma(a+b-c+1)\Gamma(c-1)}{\Gamma(a)\Gamma(b)}F(b-c+1, a-c+1, 2-c; z) \\
&=&-\frac{1}{K_2} \varphi_1(z)+z^{1-c}\frac{\Gamma(a+b-c+1)\Gamma(c-1)}{\Gamma(a)\Gamma(b)}F(b-c+1, a-c+1, 2-c; z)
\een
for $|\arg(z)|<\pi$, which means $\varphi_1(z)$ and $\varphi_2(z)$ are related, so
\be \label{zeroterm}
&&K_2\varphi_2(z)+\varphi_1(z)\nn \\
&=&\frac{-K_2 z^{1-c}}{1-c}\frac{\Gamma(a+b-c+1)\Gamma(c)}{\Gamma(a)\Gamma(b)}F(b-c+1, a-c+1, 2-c; z)
\ee
and
\be \label{limit in between}
\lim_{z\rightarrow 0}\frac{K_2\varphi_2(z)+\varphi_1(z)}{z^{1-c}}=
\frac{-K_2}{1-c} \frac{\Gamma(a+b-c+1)\Gamma(c)}{\Gamma(a)\Gamma(b)}.
\ee
Near the origin, by \eqref{geometric higher derivative}, for $n\geq 1$,
\be \label{deriv of function phi}
&&\partial^n\varphi_2(z)=\overline{\bar{\partial}^n\varphi_2(\bar{z})} \nn \\
&=&\frac{(a)_n(b)_n}{(a+b-c+1)_n}(-1)^n F(a+n,b+n,a+b-c+1+n;1-z).
\ee
By property \eqref{asym terms},
\be \label{F1}
&&F(a+n,b+n,a+b-c+1+n;1-z) \nn \\
&=&\frac{\Gamma(a+b-c+1+n)\Gamma(1-c-n)}{\Gamma(b-c+1)\Gamma(a-c+1)}F(a+n, b+n, c+n; z) \nn \\
&&+z^{1-c-n}\frac{\Gamma(a+b-c+1+n)\Gamma(c+n-1)}{\Gamma(a+n)\Gamma(b+n)}F(b-c+1, a-c+1, 2-c-n; z) \nn
\ee
for $|\arg(z)|< \pi$, then near the origin, substituting the above into \eqref{deriv of function phi}, we have
\be \label{partial k phi2}
&&\partial^n\varphi_2(z) \nn \\
&=&\frac{(a)_n(b)_n}{(a+b-c+1)_n}(-1)^n  \frac{\Gamma(a+b-c+1+n)\Gamma(1-c-n)}{\Gamma(b-c+1)\Gamma(a-c+1)}F(a+n, b+n, c+n; z) \nn \\
&&+\frac{(-1)^n}{z^{n+c-1}} \frac{\Gamma(a+b-c+1)\Gamma(c+n-1)}{\Gamma(a)\Gamma(b)} F(b-c+1, a-c+1, 2-c-n; z),
\ee
which leads to the limit
\be \label{phi2k}
\lim_{z\rightarrow 0} z^{n+c-1} \pa^n \varphi_2(z)=(-1)^n (c)_{n-1} \frac{\Gamma(a+b-c+1)\Gamma(c)}{\Gamma(a)\Gamma(b)}.
\ee
Letting $n=1$ in \eqref{partial k phi2} and combining with \eqref{zeroterm}, we have
\ben
\lim_{z\rightarrow 0}(K_2\varphi_2(\bar{z})+\varphi_1(\bar{z}))\pa \varphi_2(z)=0
\een
if $0<c=\alpha<\frac{1}{2}$. Thus
\ben
\lim_{z\rightarrow 0} \pa M(z)=\lim_{z\rightarrow 0}(K_1\varphi_1(\bar{z})+\varphi_2(\bar{z}))\pa \varphi_1(z)=\frac{ab}{c}\left(K_1-\frac{1}{K_2}\right)
\een
provided \eqref{phi11} and \eqref{former adder}.

For $(2)$, we note that \eqref{zeroterm} and \eqref{partial k phi2} is still valid for $n=1$, $\alpha=1/2$, combining with \eqref{phi11} and \eqref{phi1 and phi2} we have $(2)$ hold.

For case $(3)$,
 \be \label{partial k M}
\partial^n M(z)=(K_1\varphi_1(\bar{z})+\varphi_2(\bar{z}))\partial^n \varphi_1(z)+(K_2\varphi_2(\bar{z})
+\varphi_1(\bar{z}))\partial^n \varphi_2(z).
\ee
From properties \eqref{geometric higher derivative},
\be \label{phi1k} \displaystyle
\partial^n\varphi_1(0)=\frac{(a)_n(b)_n}{(c)_n} \ee
for $n\geq1$. Since $n>2\alpha-2$ for all $n\geq 2$ and $0<\alpha<1$, from \eqref{phi1k} and \eqref{former adder}, we know that the limit (3) is only decided by the term $(K_2\varphi_2(\bar{z})+\varphi_1(\bar{z}))\partial^n \varphi_2(z)$. Combining with \eqref{partial k M}, \eqref{phi2k} and \eqref{limit in between}, we have
\ben
&&\lim_{z\rightarrow 0}z^n |z|^{2\alpha-2}\partial^n M(z)=\lim_{z\rightarrow 0} z^n |z|^{2\alpha-2}(K_2\varphi_2(\bar{z})+\varphi_1(\bar{z})) \pa^n \varphi_2(z)\\
&=&\lim_{z\rightarrow 0}\frac{K_2\varphi_2(\bar{z})+\varphi_1(\bar{z})}{\bar{z}^{1-\alpha}}\frac{z^n }{z^{1-\alpha}} \pa^n\varphi_2(z)\\
&=&\lim_{z\rightarrow 0} \frac{K_2\varphi_2(\bar{z})+\varphi_1(\bar{z})}{\bar{z}^{1-c}} \cdot
\lim_{z\rightarrow 0} z^{n+c-1} \pa^n\varphi_2(z)\\
&=&\frac{(-1)^{n-1}(c)_{n-1}K_2}{1-c}\left(\frac{\Gamma(c)\Gamma(a+b-c+1)}{\Gamma(a)\Gamma(b)}\right)
\left(\frac{\Gamma(a+b-c+1)\Gamma(c)}{\Gamma(a)\Gamma(b)}\right)\\
&=&\frac{(-1)^{n-1}(c)_{n-1}K_2}{1-c}\left(\frac{\Gamma(c)\Gamma(a+b-c+1)}{\Gamma(a)\Gamma(b)}\right)^2
\een
as in (3).

For (4), if $m\geq1,\ n\geq 1$, we have
\ben
\bar{\pa}^m\pa^n M(z)=(K_1\bar{\pa}^m\varphi_1(\bar{z})+\bar{\pa}^m\varphi_2(\bar{z}))\pa^n \varphi_1(z)+(K_2\bar{\pa}^m\varphi_2(\bar{z})+\bar{\pa}^m\varphi_1(\bar{z}))\pa^n \varphi_2(z).
\een
Since
\ben
\lim_{z \rightarrow 0}z^{n+c-1} \pa^{n}\varphi_1(z)=0,
\een
then
\ben
&&\lim_{z\rightarrow 0}\bar{z}^m z^n |z|^{2\alpha-2}\bar{\pa}^m \pa^n M(z)=\lim_{z\rightarrow 0}\frac{\bar{z}^m z^n}{|z|^{2-2c}}K_2\bar{\pa}^m\varphi_2(\bar{z})\pa^n \varphi_2(z)\\
&=&\lim_{z\rightarrow 0}K_2 \frac{\bar{z}^m}{\bar{z}^{1-c}}\bar{\pa}^m\varphi_2(\bar{z})\cdot\frac{z^n}{z^{1-c}}\pa^n \varphi_2(z)\\
&=&(-1)^{m+n}(c)_{m-1}(c)_{n-1}K_2\left(\frac{\Gamma(c)\Gamma(a+b-c+1)}{\Gamma(a)\Gamma(b)}\right)^2
\een
as in (4). \hfill $\Box$\\

The following result is a specific version of Theorem \ref{general pu}.
\begin{theorem} \label{pu}
\textsl{For $\lambda(z):=\lambda_{\alpha,\,\beta,\,\gamma}(z)$ as in \eqref{lambda(z)1} with order $\alpha\in(0,\,1)$, let $u(z):=\log\lambda(z)$. Then for $m,\,n\geq 1$, \vspace*{2mm} \\
(i) $ \displaystyle \ \lim_{z\rightarrow0}z^n\partial^nu(z)=\frac {\alpha}2(-1)^n(n-1)!=\lim_{z\rightarrow0}\bar{z}^n\bar{\pa}^nu(z),$ \vspace*{2mm} \\
(ii) $ \displaystyle \ \lim_{z\rightarrow 0}\bar{z}^m z^n |z|^{2\alpha-2}\bar{\pa}^m\pa^nu(z)=\frac{(-1)^{n+m}(c)_{n-1}(c)_{m-1}K^2_2}{K_1K_2-1} \left(\frac{\Gamma(c)\Gamma(a+b-c+1)}{\Gamma(a)\Gamma(b)}\right)^2.$}
\end{theorem}
\textbf{Remark.} Theorem \ref{general pu} was proved for the order $0<\alpha \leq 1$ in \cite{Zhang1} with a different limit for the mixed differential, while Theorem \ref{pu} is given for the order $0<\alpha<1$ and we prove it in a different way for the completeness of this paper. The proof of Theorem \ref{pu} also can be taken to be an application of Lemma \ref{Mz lemma}. For the hyperbolic density $\lambda(z)$ with order $\alpha=1$, we can also prove Theorem \ref{general pu} directly by discussing the properties of hypergeometric functions. \vspace*{3mm} \\
\textbf{Proof of Theorem \ref{pu}.} We note that $$u(z)=-\alpha\log|z|-\beta\log|1-z|+\log K_3-\log M(z)$$ with $M(z)$ as in \eqref{Mz}. At first we consider $\partial^n\log M(z)$. From \eqref{phi1 and phi2},
\be \label{m0}
M(0)=K_1\varphi_1(0)+\varphi_2(0)=K_1-\frac 1 {K_2}.
\ee
We can calculate that $M(0)>0$ for any $a,\ b,\ c$ as in \eqref{abc}. From Lemma \ref{Mz lemma} we have
$$\displaystyle  \lim_{z\rightarrow0}z^k\partial^kM(z)=0$$
for all $k\geq1$ and $0<\alpha<1$.
It is easy to observe that $\partial^n\log M(z)$ is a linear combination of products of  $\displaystyle \frac {\partial^kM}{M}$ with $k\leq n$, so when $n \geq 1$, $$\lim_{z\rightarrow0}z^n\partial^n\log M(z)=0.$$
Since $$\displaystyle \partial^n\log|1-z|=-\frac{(n-1)!}{2(1-z)^n},\ \ \
\partial^n\log|z|=\frac{(-1)^{n-1}(n-1)!}{2z^n},$$ then the first equality in (i) holds.

For the second equality, note that $u(z)$ is real-valued,
$$\displaystyle \ \lim_{z\rightarrow0}\bar{z}^n\bar{\partial}^nu(z)
=\lim_{z\rightarrow0}\overline{z^n\partial^nu(z)}=\frac {\alpha}2(-1)^n(n-1)!.$$
Therefore (i) is valid.

Now we discuss the term $\bar{\partial}^m\partial^n\log M(z)$ to complete the proof. Since $\bar{\partial}^m\partial^n\log M(z)$ is a linear combination of products of $\bar{\pa}^t \pa^k M / M$ with $0\leq t\leq m$, $0\leq k\leq n$, so Lemma \ref{Mz lemma} implies that
$$\lim_{z\rightarrow0}\bar{z}^mz^n|z|^{2\alpha-2}\prod^{N}_{j=2}\frac{\bar{\pa}^{t_j}\pa^{k_j}M(z)}{M(z)}=0,$$
where $2 \leq N \leq m+n$, $1 \leq t_j \leq m$ and $1 \leq k_j \leq n$ for every index $j$, $2 \leq j \leq N$. Thus
\ben
&&\lim_{z\rightarrow 0}\bar{z}^m z^n |z|^{2\alpha-2}\bar{\pa}^m\pa^n\log M(z)=\lim_{z\rightarrow 0}\bar{z}^m z^n |z|^{2\alpha-2}\frac{\bar{\pa}^m\pa^n M(z)}{M(z)} \\
&=&\frac{(-1)^{n+m}(c)_{n-1}(c)_{m-1}K^2_2}{K_1K_2-1} \left(\frac{\Gamma(c)\Gamma(a+b-c+1)}{\Gamma(a)\Gamma(b)}\right)^2.
\een
Note that $\bar{\partial}^m\partial^n\log|1-z|=0$, $\bar{\partial}^m\partial^n\log|z|=0$, thus (ii) holds. \hfill $\Box$
\vspace{2mm}
\par For the hyperbolic metric, the following result corresponding to Theorem \ref{coro}  holds.
\begin{theorem} \label{lc}
\textsl{For $m,\,n \geq 0$, $0<\alpha<1$ and $\lambda(z)$ as in (\ref{lambda(z)1}), the limit
$$l_{m,\,n}:=\frac{1}{m!n!}\lim_{z \rightarrow0}|z|^{\alpha}{\bar{z}}^m z^n {\bar{\partial}}^m \partial^n \lambda(z)$$
exists. Let
\be \label{l00}
l_{0,\,0}=l:=\lim_{z \rightarrow0}|z|^{\alpha}\lambda(z)=\frac{\delta}{1-\delta^2}(1-\alpha)
\ee
by Theorem \ref{limit for alpha less than 1}, then the numbers $l_{m,\,n}$ satisfy the following  \vspace*{1mm} \\
(i) \,$l_{m,\,n}=\displaystyle{{-\frac{\alpha}{2} \choose n}{-\frac{\alpha}{2} \choose m}l}$, \vspace*{1mm} \\
(ii) $l_{m,\,n}=l_{n,\,m}$,\\
where
$${\tau \choose j}=\frac{\tau(\tau-1)\cdots (\tau-j+1)}{j\;!}$$
is the binomial coefficient. }
\end{theorem}
\textbf{Proof.}
Since
\be \label{pa1 lam}
\pa \lambda (z)=\lambda(z)\, \pa u(z)
\ee
we have
$$\pa^n\lambda(z)=\sum_{j=0}^{n-1}{n-1 \choose j}\pa^{n-j}u(z)
\,\pa^j\lambda(z)$$
by induction, where $\pa^0\lambda(z)=\bar{\pa}^0\lambda(z)=\lambda(z).$
Then$$l_{0,\,n}=\frac{1}{n!}\lim_{z\rightarrow0}\sum_{j=0}^{n-1}{n-1 \choose j}
z^{n-j}\pa^{n-j}u(z)\cdot|z|^{\alpha}z^j\pa^j\lambda(z).$$
From the existence of $\lim_{z\rightarrow0}z^{n-j}\pa^{n-j}u(z)$ and $l$, it is known that $l_{0,\,n}$ exists. By (ii) in Theorem \ref{pu}, we have $$\lim_{z\rightarrow 0}\bar{z}^m z^n\bar{\pa}^m \pa^n u(z)=0.$$ So we can write $l_{m,\,n}$ as a sum of the terms not containing any mixed derivatives of $u(z)$,
\be\label{ind}
l_{m,\,n}=\frac{1}{m!n!}\lim_{z\rightarrow0}\sum_{j=0}^{n-1}{n-1 \choose j}
z^{n-j}\pa^{n-j}u(z)\,|z|^{\alpha}\bar{z}^m z^j\bar{\pa}^m\pa^j\lambda(z),
\ee
thus the existence of $l_{0,\,n}$ guarantees $l_{m,\,n}$ exists.

If $m=0$, $n=1$, then \eqref{l00} and \eqref{pa1 lam} give
$$l_{0,1}=\lim_{z\rightarrow 0}|z|^{\alpha}z\pa\lambda(z)=\lim_{z\rightarrow 0}|z|^{\alpha}\lambda(z) \cdot z\pa u(z)=-\frac{\alpha}{2}l,$$
which is a real number, so $l_{1,\,0}=\overline{l_{0,\,1}}=l_{0,\,1}$.
Note that
\be \label{pa bar n lam}
{\displaystyle\bar{\pa}^n\lambda(z)=\sum_{j=0}^{n-1}{n-1 \choose j}
\bar{\pa}^{n-j}u(z)\,\bar{\pa}^j\lambda(z)},
\ee
then $l_{n,\,0}=l_{0,\,n}$ by induction. From \eqref{ind}, \eqref{pa bar n lam}, and (i) of Theorem \ref{pu}, we have
\ben
l_{m,\,n}&=&\sum_{j=0}^{n-1} \lim_{z\rightarrow0} \frac{1}{m!n!} \frac{(n-1)!}{j!(n-1-j)!}
z^{n-j}\pa^{n-j}u(z)\cdot|z|^{\alpha}\bar{z}^m z^j\bar{\pa}^m\pa^j\lambda(z) \\
&=&\frac{1}{n}\sum_{j=0}^{n-1}  \frac{1}{m!} \frac{1}{j!(n-1-j)!}\lim_{z\rightarrow0}
z^{n-j}\pa^{n-j}u(z)\cdot\lim_{z\rightarrow0}|z|^{\alpha}\bar{z}^mz^j\bar{\pa}^m\pa^j\lambda(z) \\
&=&\frac{1}{n}\sum_{j=0}^{n-1} \frac{\alpha(-1)^{n-j}}{2}\frac{1}{m!j!} \lim_{z\rightarrow0}|z|^{\alpha}\bar{z}^m z^j\bar{\pa}^m\pa^j\lambda(z)
=\frac{\alpha}{2n}\sum_{j=1}^{n-1}(-1)^{n-j}l_{m,\,j}.
\een
Then
$$ n \cdot l_{m,\,n}=\frac{\alpha}2\sum_{j=0}^{n-2}(-1)^{n-j}l_{m,\,j}-\frac{\alpha}{2}l_{m,\,n-1}
=-(n-1)l_{m,\,n-1}-\frac{\alpha}{2}l_{m,\,n-1}.$$
Since $l_{0,\,n}=l_{n,\,0}$,
\begin{eqnarray*}
l_{m,\,n}&=&\frac{-\frac{\alpha}{2}-n+1}{n}l_{m,\,n-1}=
{-\frac{\alpha}{2} \choose n}l_{m,\,0}\\
&=&{-\frac{\alpha}{2} \choose n}l_{0,\,m}={-\frac{\alpha}{2} \choose n}{-\frac{\alpha}{2} \choose m}l_{0,\,0}.
\end{eqnarray*}
Thus (i) holds and (ii) follows immediately form (i). \hfill $\Box$\\

The following estimate is for the general case.
\begin{theorema} $(\mathrm{[15]})$ \label{estimate v}
\textsl{ Let $\kappa:\mathbb{D}\rightarrow\mathbb{R}$ be a locally H\"{o}lder continuous function with $\kappa(0)<0$. If $u:\mathbb{D}^*\rightarrow \mathbb{R}$ is a $C^2$-solution to $\Delta u=-\kz e^{2u}$ in $\mathbb{D}^*$, then $u$ has an order $\alpha \in (-\infty,\,1]$. If, in addition, $\kz \in C^{n-2,\,\nu}(\mathbb{D}^*)$ for an integer $n \geq 3$, $0<\nu \leq 1$, then $u(z) \in C^{n,\,\nu}(\mathbb{D}^*)$ by the regularity theorem. If the order $0<\alpha<1$, then for the remainder function $v(z)$ and for $n_1,\,n_2 \geq 1$, $n_1+n_2 =n$, near the origin, we have
$$\partial^n v(z), \ \bar{\partial}^n  v(z), \ \bar{\partial}^{n_1}\partial^{n_2}v(z)=\textit{O}(|z|^{2-2\alpha-n}).$$}
\end{theorema}

\par From the proof of Theorem \ref{pu}, we can provide a way to verify the sharpness of Theorem \ref{estimate v}, and also Theorem 1.1 in \cite{Rothbehaviour}. We state the result as the following theorem.

\begin{theorem} \label{vn}
\textsl{For $m,\,n \geq 1$ and $\lambda$ as in (\ref{lambda(z)1}) with the order $0<\alpha<1$, then near the origin, the remainder function $v(z)$ satisfies \vspace*{2mm}\\
(i)$\displaystyle \lim_{z\rightarrow 0} \pa v(z)=\frac{ab}{c}$ for $\displaystyle \ 0<\alpha<{1}/{2}$, \vspace*{2mm}\\
(ii)$\displaystyle \pa v(z)=2ab+\frac{2K^2_2}{K_1K_2-1} \left(\frac{\Gamma(c)\Gamma(a+b-c+1)}{\Gamma(a)\Gamma(b)}\right)^2 \frac{\bar{z}}{|z|} +O\left(|z|^{\frac{1}{2}}\right)$ \vspace{1mm}\\near the origin for $ \alpha={1}/{2}$, \vspace*{2mm}\\
(iii)$\displaystyle \lim_{z\rightarrow 0} z^{n}|z|^{2\alpha-2}\pa^{n}v(z)=\frac{(-1)^{n-1}(c)_{n-1}K_2^2}{(1-c)(K_1K_2-1)} \left(\frac{\Gamma(c)\Gamma(a+b-c+1)}{\Gamma(a)\Gamma(b)}\right)^2$ \vspace{1mm}\\for $n\geq 2$ if $\ 0<\alpha\leq {1}/{2}$ and $n\geq 1$ if $\ {1}/{2}<\alpha<1$, \vspace*{2mm}\\
(iv)$\displaystyle \lim_{z\rightarrow 0} \bar{z}^{m}z^{n}|z|^{2\alpha-2}\bar{\pa}^{m}\pa^{n}v(z)=\frac{(-1)^{n+m}(c)_{n-1}(c)_{m-1}K^2_2}{K_1K_2-1} \left(\frac{\Gamma(c)\Gamma(a+b-c+1)}{\Gamma(a)\Gamma(b)}\right)^2$ \vspace{1mm}\\for $m,\,n \geq 1$ and $0<\alpha<1$.}
\end{theorem}
\textbf{Proof.} Since for $\lambda(z)$ in (\ref{lambda(z)1}), $v(z)=-\beta \log|1-z|+ \log K_3 - \log M(z)$ and $$\partial^n\log|1-z|=\frac{-(n-1)!}{2(1-z)^n},$$ we consider $\partial^n \log M(z)$ only. From the proof of Theorem \ref{pu}, the limits $ \lim_{z\rightarrow 0} \pa v(z)$ and $\lim_{z\rightarrow 0}z^n |z|^{2\alpha-2}\pa^n v(z)$ both depend solely on the term $\partial^n M(z)$.
Thus by Lemma \ref{Mz lemma} and \eqref{m0},
$$ \lim_{z\rightarrow 0} \pa v(z)=\lim_{z\rightarrow 0} \frac{\pa M(z)}{M(z)},$$
$$\lim_{z\rightarrow 0}z^n |z|^{2\alpha-2}\partial^n v(z)=\lim_{z\rightarrow 0}z^n |z|^{2\alpha-2}\frac{\partial^n M(z)}{M(z)}.$$
So we obtain the four cases above corresponding to ones in Lemma \ref{Mz}. \hfill $\Box$

\section{Case $\alpha =1$}
\setcounter{equation}{0}

If $\alpha =1$, the formula for $\lambda_{1,\,\beta,\gamma}$ is to be understood in the limit sense $\lim_{\alpha \rightarrow 1{-}}$. So when $\alpha= c=1$, we have
\be
K_3=\frac{1}{B( a, b)}:=\frac 1 {B}, \ \ K_2=0, \nn \\
S:=\frac{\pi\sin(\pi( a+ b))}{\sin\pi a\sin\pi b}, \ \ K_1=-\frac S B, \label{k1}
\ee
$$\varphi_1(z)=F( a, b,1;z),\ \ \ \varphi_2(z)=F( a, b, a+ b;1-z).$$
Then
\be
\lambda_{1,\,\beta,\gamma}(z)&=&\frac 1 {|z|}\frac 1 {|1-z|^{\beta}} \frac{K_3}{K_1|\varphi_1(z)|^2+\varphi_1(z)\varphi_2(\bar{z})
+\varphi_1(\bar{z})\varphi_2(z)} \label{lambda for alpha=1}\\
& :=&\frac 1 {|z|} \frac 1 {|1-z|^{\beta}} \frac {K_3}{M(z)}, \nn
\ee
and the remainder function of $u(z)$ near the origin is
\be \label{remainder w}
w(z)=-\beta\log|1-z|+\log {K_3}-\log M(z)+\log\log(1/|z|).
\ee
The assumption of Theorem \ref{explicit formula} and \eqref{abc} show that $a$ and $b$ satisfy
$$0<a<1,\quad 0<b<1/2,\quad 0<a+b<1.$$
The function
\be \label{2R-S}
2R-S=4\Psi(1)-2\Psi( a)-2\Psi( b)-\pi \cot\pi a-\pi \cot\pi b
\ee
is of special interest where $R:=R(a,b)$ is as in \eqref{rab} and $S$ is given by \eqref{k1}.
Let $G(x):=2\left(\Psi(1)-\Psi(x)\right)- \pi \cot\pi x$. For the Gamma function $\Gamma$ and $0<x<1$, we have
$\Gamma(x)\Gamma(1-x)=\pi / \sin \pi x$. Taking the logarithmic derivatives of both sides leads to
$$\frac{\Gamma'(x)}{\Gamma(x)}-\frac{\Gamma'(1-x)}{\Gamma(1-x)}=-\pi \cot \pi x.$$
So
\be \label{Gx}
G(x)=2\Psi(1)-\Psi(x)-\Psi(1-x),
\ee
which means $G(x)=G(1-x)$. The fact that the digamma function is negative and decreasing on $(0,\,1)$ implies that $G(x)>0$ when $0<x<1$. Since $2R-S=G(a)+G(b)$, then $2R-S>0$ for all $a,\ b$ given by \eqref{abc}. \vspace*{2mm}

For any conformal metric $\lambda(z)|dz|$ with the negative Gaussian curvature and the remainder function $w(z)$ defined by \eqref{singularity}, we have the following result to describe the asymptotic behavior of $w(z)$ near the origin.
\begin{theorema} $(\mathrm{[15]})$ \label{generalized hyper w}
\textsl{ Let $\kappa:\mathbb{D}\rightarrow\mathbb{R}$ be a locally H\"{o}lder continuous function with $\kappa(0)<0$. If $u:\mathbb{D}^*\rightarrow \mathbb{R}$ is a $C^2$-solution to $\Delta u=-\kz e^{2u}$ in $\mathbb{D}^*$, then $u$ has an order $\alpha \in (-\infty, 1]$. If, in addition, $\kz \in C^{n-2,\,\nu}(\mathbb{D}^*)$ for an integer $n \geq 3$, $0<\nu \leq 1$, then $u(z) \in C^{n,\,\nu}(\mathbb{D}^*)$ by the regularity theorem. If order $\alpha=1$, then for the remainder function $w(z)$ and for $n_1,\,n_2 \geq 1$, $n_1+n_2 =n$, near the origin, we have
\ben
\bar{\partial}^nw(z), \ \partial^nw(z)=\textit {O}(|z|^{-n}\log^{-2}(1/|z|)),  \\
\bar{\partial}^{n_1}\partial^{n_2}w(z)=\textit {O}(|z|^{-n}\log^{-3}(1/|z|)).
\een}
\end{theorema}
We can verify the sharpness of Theorem \ref{generalized hyper w} by use of $\lambda_{1,\,\beta,\,\gamma}(z)$ as in \eqref{lambda for alpha=1}. Furthermore, for $\lambda(z)$, we can obtain its precise estimate for higher order derivatives of $w(z)$ near the origin. In fact, we have the following result stronger than Theorem \ref{generalized hyper w}.
\begin{theorem} \label{estiofw}
\textsl{Let $\lambda(z):=\lambda_{1,\,\beta,\,\gamma}(z)$ as in \eqref{lambda for alpha=1} with $\beta$ and $\gamma$ satisfying the condition in Theorem \ref{explicit formula}, and $w(z)$ be the remainder function as in \eqref{remainder w}. Then for $m,\: n\geq1$, we have \vspace{2mm}\\
(i)$\displaystyle \ \lim_{z\rightarrow 0}z^{n}\log^2(1/|z|)\partial^nw(z)=\frac {(-1)^n(n-1)!}{4}(G(a)+G(b)),$ \vspace{2mm}\\
(ii)$ \displaystyle \  \lim_{z\rightarrow 0}z^n\bar{z}^m\log^3(1/|z|)\bar{\partial}^m\partial^nw(z)=\frac {(-1)^{m+n-1}(n-1)!(m-1)!}{4}(G(a)+G(b)),$ \vspace{2mm} \\
where the function $G$ is defined by \eqref{Gx} and $a,\:b$ are given by \eqref{abc}.}
\end{theorem}
\textbf{Proof.} \ For the remainder function given by \eqref{remainder w}, we discuss $\log\log(1/|z|)$ and $\log M(z)$ separately. At first, consider the higher order derivatives of $\log\log(1/|z|)$. By induction we know that
$$\partial^n\log\log(1/|z|)=\sum^n_{j=1}\frac{C^{(n)}_j}{z^n\log^j(1/|z|)}$$
with constant $C^{(n)}_{j}$ for $1 \leq j \leq n$.
Here we only need the first two terms of $\pa^{n}\log\log(1/|z|)$ for future use. As for the pure derivative  $\partial^n\log\log(1/|z|)$ with $n\geq 1$, set $\mathcal{A}_n:=C^{(n)}_1$ and $\mathcal{B}_n:=C^{(n)}_2$, so
$$\partial^n\log\log(1/|z|)=\frac{\mathcal{A}_n}{z^n\log(1/|z|)}
+\frac{\mathcal{B}_n}{z^n\log^2(1/|z|)}+\sum^n_{j=3}\frac{C^{(n)}_j}{z^n\log^j(1/|z|)},$$
then the following recurrent relations hold,
$$\mathcal{A}_1=-\frac 1 2, \ \ \  \mathcal{B}_1=0,$$
$$\mathcal{A}_n=-(n-1)\mathcal{A}_{n-1},\ \ \ \mathcal{B}_n=-(n-1)\mathcal{B}_{n-1}+\frac 1 2 \mathcal{A}_{n-1}.$$
Thus
\be
\mathcal{A}_n=\frac{(-1)^n}{2}(n-1)!, \qquad \label{an} \\
\mathcal{B}_n=\frac{(-1)^{n-1}}{4}(n-1)!\sum^{n-1}_{j=1}\frac 1 j.  \label{bn}
\ee
For the mixed derivative case with $n\geq1$, $m\geq1$, we fix $n$, so by induction,
$$\bar{\pa}^{m}\pa^n\log\log(1/|z|)=\sum^m_{j=1}\frac{C^{(m,\,n)}_{j}}{\bar{z}^m z^n\log^{j+1}(1/|z|)}$$
with constant $C^{(m,\,n)}_{j}$ for $1 \leq j \leq m$. Set $\mathcal{C}_m:=C^{(m,\,n)}_1$ and $\mathcal{D}_m:=C^{(m,\,n)}_2$, we have
$$\bar{\partial}^m \partial^n\log\log(1/|z|)=\frac{\mathcal{C}_m}{\bar{z}^m z^n\log^2(1/|z|)}
+\frac{\mathcal{D}_m}{\bar{z}^m z^n\log^3(1/|z|)}+\sum^m_{j=3}\frac{C^{(m,\,n)}_{j}}{\bar{z}^m z^n\log^{j+1}(1/|z|)}.$$
Then
$$\mathcal{C}_1=\frac 1 2 \mathcal{A}_n, \ \ \  \mathcal{D}_1=\mathcal{B}_n,$$
$$\mathcal{C}_m=-(m-1)\mathcal{C}_{m-1},\ \ \ \mathcal{D}_m=-(m-1)\mathcal{D}_{m-1}+\mathcal{C}_{m-1}.$$
Therefore
\be
\mathcal{C}_m=\frac{(-1)^{m+n-1}}{4}(m-1)!(n-1)!, \qquad \qquad \label{cn} \\
\mathcal{D}_m=\frac{(-1)^{m+n}}{4}(m-1)!(n-1)!\left(\sum^{n-1}_{j=1}\frac 1 j +\sum^{m-1}_{j=1}\frac 1 j \right). \label{dn}
\ee

To estimate the derivatives of $\log M(z)$, we first calculate $\displaystyle \frac{\partial^n M(z)}{M(z)}$ for $n\geq 1$. Since
\be\begin{array}{l}
\displaystyle\partial^n\varphi_1(z)=\frac{(a)_n(b)_n}{n!}F(a+n,b+n,n+1;z),\\
\displaystyle\partial^n\varphi_2(z)=(-1)^n
\frac{(a)_n(b)_n}{(a+b)_n}F(a+n,b+n,a+b+n;1-z),\nonumber
\end{array}\ee
so $$\partial^n\varphi_1(0)=\frac{(a)_n(b)_n}{n!},$$
and near the origin, by \eqref{behavior of hypergeometric} we have
$$\partial^n\varphi_2(z)= \frac{(a)_n(b)_n}{(a+b)_n}\frac{(-1)^n}{z^n}F(b,a,a+b+n;1-z).$$
Considering \eqref{covering case} shows that
$$F(b,a,a+b+n;1-z)=\frac{\Gamma(a+b+n)\Gamma(n)}{\Gamma(a+n)\Gamma(b+n)}+\textit{O}(|z|\log|z|)$$
near the origin. Thus
\be
\partial^n\varphi_2(z)&=&\frac{(-1)^n}{z^n} \frac{(a)_n(b)_n}{(a+b)_n}
\left(\frac{\Gamma(a+b+n)\Gamma(n)}{\Gamma(a+n)\Gamma(b+n)}+\textit{O}(|z|\log|z|)\right) \nn \\
&=&\frac{(-1)^n (n-1)!}{B z^n}+\textit{O}\left(\frac{\log|z|}{|z|^{n-1}}\right). \nn
\ee
Property \eqref{behavior of hypergeometric} gives
$$\varphi_2(z)=\frac{1}{B}\left(\log \frac{1}{z}+R\right)\left(1+O(z)\right),$$
summary the estimates above, so we can obtain
\ben
\pa^nM(z)&=&(K_1\varphi_1(\bar{z})+\varphi_2(\bar{z}))\pa^n \varphi_1(z)+\varphi_1(\bar{z})\pa^n \varphi_2(z) \\
&=&\frac{(a)_n(b)_n}{B\ n!}\log\frac{1}{\bar{z}}+\frac{(-1)^n (n-1)!}{B z^n}+\textit{O}\left(\frac{\log|z|}{|z|^{n-1}}\right)
\een
and
\ben
M(z)&=&K_1\varphi_1(\bar{z})\varphi_1(z)+\varphi_2(\bar{z})\varphi_1(z)+\varphi_1(\bar{z}) \varphi_2(z)\\
&=&\frac{2\log(1/|z|)}{B} \left(1+\frac{2R-S}{2\log(1/|z|)}+\textit {O}(|z|)\right)
\een
near the origin. Then
\be \label{mz}
&&\frac{\partial^nM(z)}{M(z)} \nn \\
&=&\pa^n M(z)\ \frac{B}{2\log(1/|z|)}\left(1-\frac{2R-S}{2\log(1/|z|)}+\textit {O}(|z|)\right) \nn \\
&=&\frac{(-1)^{n}(n-1)!}{2z^n\log(1/|z|)}-\frac{(-1)^{n}(n-1)!(G(a)+G(b))}
{4z^n\log^2(1/|z|)}+\textit{O}\left(\frac 1 {|z|^{n-1}}\right). \quad \quad
\ee
We note that
\ben
&&\bar{\pa}^m \pa^nM(z) \\
&=&(K_1\bar{\pa}^m\varphi_1(\bar{z})+\bar{\pa}^m \varphi_2(\bar{z}))\pa^n \varphi_1(z)+\bar{\pa}^m \varphi_1(\bar{z})\pa^n \varphi_2(z) \\
&=&\frac{(-1)^m (m-1)!}{B \bar{z}^m}\frac{(a)_n(b)_n}{n!}+\frac{(-1)^n (n-1)!}{B z^n}\frac{(a)_m(b)_m}{m!}+\textit{O}\left(\frac{\log|z|}{|z|^{m-1}}\right)+\textit{O}\left(\frac{\log|z|}{|z|^{n-1}}\right),
\een
the same technique leads to
\be \label{mixed}
\frac{\bar{\partial}^m\partial^n M(z)}{M(z)}=\textit{O}\left(\frac 1 {|z|^{\tau}\log(1/|z|)}\right),
\ee
where $\tau=\max\{m,\,n\}< m+n.$

\vspace{1mm}
\par Now we can consider derivatives of $w(z)$. In the pure derivative case,
\be \label{partial n w}
\partial^nw(z)=\frac {\beta(n-1)!}{2(1-z)^n}-\partial^n\log M(z)+\partial^n\log\log(1/|z|).\ee
Since the coefficients of the first two differential terms in $\partial^n\log\log(1/|z|)$ are already known as $\mathcal{A}_n$ and $\mathcal{B}_n$, now we discuss $\partial^{n}\log M(z)$. Note that $\partial^{n}\log M(z)$ is a linear combination of finitely many terms of the form
\be \label{oneterm}
\prod^{k}_{j=1}\frac{\pa^{n_j}M(z)}{M(z)}
\ee
for $1\leq k \leq n$. When $k=1$, term \eqref{oneterm} is corresponding to the first term in the third line of \eqref{mz} with $n=1$, and it will be canceled by $\mathcal{A}_n$ given in \eqref{an}. So we should look at the second term which contains $z^{-n}\log^{-2}(1/|z|)$ for the higher order derivatives, while the higher power terms in \eqref{mz} are ignored for a moment. For \eqref{oneterm}, estimate \eqref{mz} shows that
$$\prod^{k}_{j=1}\frac{\pa^{n_j}M(z)}{M(z)}=\textit{O}\left(\frac{1}{z^n \log^k (1/|z|)}\right)\ {\textrm{for}} \ n=\sum^{k}_{j=1} n_j,$$
therefore, to generate the $z^{-n}\log^{-2}(1/|z|)$ term, $k$ is at most 2, thus the $z^{-n}\log^{-2}|z|$ term of $\partial^n \log M(z)$ only appears in
$$\frac{\partial^n M}{M}-\frac 1 2 \sum_{j=1}^{n-1}{n \choose j} \frac{\partial^jM\partial^{n-j}M}{M^2}.$$
For every $1\leq j\leq n-1$, we have
$$\frac{\partial^jM\partial^{n-j}M}{M^2}=\frac{(-1)^n(j-1)!(n-j-1)!}{4z^n\log^2|z|}+
\textit {O}(\frac{1}{|z|^{n}\log^{3}(1/|z|)}).$$
Denote the coefficient of $ z^{-n}\log^{-2}|z|$ in $\sum_{j=1}^{n-1}{n \choose j} (\pa^j M\partial^{n-j}M/{M^2})$ by $k_n$, then \eqref{mz} leads to
\be\label{kn}
k_n &=&\sum^{n-1}_{j=1}{n \choose j} \frac{(-1)^n(j-1)!(n-j-1)!}{4}=\frac{(-1)^{n}(n-1)!}{4}\sum^{n-1}_{j=1}\frac{n}{j(n-j)} \nn \\
&=&\frac{(-1)^{n}(n-1)!}{4}\sum^{n-1}_{j=1}\left(\frac{1}{j}+\frac{1}{n-j}\right) \nn \\
&=&\frac{(-1)^{n}(n-1)!}{2}\sum_{j=1}^{n-1}\frac{1}{j}=-2\mathcal{B}_n.
\ee
Note that the term $z^{-n}\log^{-1}(1/|z|)$ only appears in ${\partial^n M(z)}/{M(z)}$, and (\ref{an}), (\ref{mz}), \eqref{partial n w} show that $z^{-n}\log^{-1}(1/|z|)$ actually is canceled in $\partial^n w(z)$. In combination with \eqref{bn} and \eqref{kn}, we obtain
\ben
\lim_{z\rightarrow 0}z^{n}\log^2(1/|z|)\partial^nw(z)
&=&\frac{(-1)^{n}(n-1)!(G(a)+G(b))}{4}+\frac{1}{2}k_n+\mathcal{B}_n \\
&=&\frac {(-1)^n(n-1)!}{4}(G(a)+G(b)),
\een
thus (i) holds.

For the mixed derivatives case,
\ben
\bar{\partial}^m \partial^n w(z)=-\bar{\partial}^m \partial^n\log M(z)+\bar{\partial}^m \partial^n\log\log(1/|z|).
\een
Since the coefficients of the first two terms in $\bar{\partial}^m\partial^n\log\log(1/|z|)$ are given as $\mathcal{C}_m$ and $\mathcal{D}_m$, now we consider $\bar{\partial}^m\partial^n\log M(z)$. It is known that $M(z)=M(\bar{z})$ and $\bar{\partial}^m\partial^n\log M(z)=\overline{\bar{\partial}^n\partial^m\log M(z)}$. Thus without loss of generality we may assume $m\leq n$. Similarly as in the pure derivative case, there will be some cancelation for the term containing $z^{-n}\bar{z}^{-m}\log^{-2}(1/|z|)$, so the coefficient of $z^{-n}\bar{z}^{-m}\log^{-3}(1/|z|)$ is desired. The term containing $z^{-n}\bar{z}^{-m}\log^{-3}(1/|z|)$ must the product of at most three terms in the forms of ${\partial^{n_j}M}/{M}$ or ${\bar{\partial}^{m_j}M}/{M}$. Estimate (\ref{mixed}) and \eqref{mz} imply that the term $z^{-n}\bar{z}^{-m}\log^{-3}(1/|z|)$ of $\bar{\partial}^m\partial^n\log M(z)$ only appears in
$$-\frac{\bar{\partial}^m M \partial^n M}{M^2}+\frac{\bar{\partial}^m M}{M}\sum^{n-1}_{j=1}{n \choose k} \frac{\partial^j M \partial^{n-j} M}{M^2}+\frac{\partial^n M}{M}\sum^{m-1}_{j=1}{m \choose k} \frac{\bar{\partial}^j M \bar{\partial}^{m-j} M}{M^2}$$
for $m\geq 1$, $n\geq 1$. Then
$$\bar{\partial}^m\partial^n\log M(z)=\frac {t_{m,\,n}}{z^n\bar{z}^m\log^3(1/|z|)} +\textit{O}\left(\frac{1}{|z|^{m+n+1}\log^{3}(1/|z|)}\right)$$
with
\be \label{tmn}
t_{m,\,n}=\left(G(a)+G(b)+\sum^{n-1}_{j=1}\frac{1}{j}+\sum^{m-1}_{j=1}\frac{1}{j}\right)\frac {(-1)^{m+n}(m-1)!(n-1)!}{4},
\ee
which can be obtained by the technique similar to the one applied to $k_n$ in \eqref{kn}. Note that the term $\bar{z}^{-m}z^{-n}\log^{-2}(1/|z|)$ only occurs in ${\bar{\partial}^{m}M \partial^n M}/{M^2}$ with the coefficient ${(-1)^{m+n}(n-1)!(m-1)!}/{4}$, comparing with (\ref{cn}) shows that there is no term of $\bar{z}^{-m}z^{-n}\log^{-2}(1/|z|)$ left in the expression for $\bar{\partial}^m\partial^nw(z)$.
Thus for $m\geq1,\ n\geq 1$, by (\ref{dn}) and \eqref{tmn}, we obtain
\ben
\lim_{z\rightarrow 0}z^n\bar{z}^m\log^3(1/|z|)\bar{\partial}^m\partial^nw(z)=-t_{m,\,n}+\mathcal{D}_m
=\frac {(-1)^{m+n-1}(n-1)!(m-1)!}{4}(G(a)+G(b))
\een
This completes the proof and verifies the sharpness of Theorem \ref{generalized hyper w}. \hfill $\Box$ \\

\par When the order $\alpha=1$, there is an analogue of Theorem \ref{pu} and Theorem \ref{lc}, see \cite{Zhang1}, also \cite{Rothbehaviour}. Here we list them as following ones without proof.

\begin{theorem} \label{pu2}
\textsl{For $\lambda(z):=\lambda_{1,\,\beta,\gamma}(z)$ as in \eqref{lambda for alpha=1}, let $u(z):=\log\lambda(z)$. Then for $m,\,n\geq 1$, \vspace*{2mm} \\
(i) $ \displaystyle \ \lim_{z\rightarrow0}z^n\partial^nu(z)=\frac {1}2(-1)^n(n-1)!=\lim_{z\rightarrow0}\bar{z}^n\bar{\pa}^nu(z),$ \vspace*{2mm} \\
(ii) $ \displaystyle \ \lim_{z\rightarrow 0}\bar{z}^m z^n \log^2(1/|z|)\bar{\pa}^m\pa^nu(z)=\frac{(-1)^{n+m}(n-1)!(m-1)!}{4}.$}
\end{theorem}

\begin{theorema} $(\mathrm{[15]})$ \label{lctilte}
\textsl{For $m,\,n \geq 0$, $\alpha=1$ and $\lambda(z)$ as in \eqref{lambda for alpha=1}, the limit
$$l'_{m,\,n}:=\frac{1}{n!m!}\lim_{z \rightarrow0}|z|\log(1/|z|){\bar{z}}^m z^n {\bar{\partial}}^m \partial^n \lambda(z)$$
exists. Moreover, the numbers $l'_{m,\,n}$ satisfy the following \vspace*{1mm}\\
(i)$ \:\:\:\,\displaystyle  l'_{0,\,0}:=\lim_{z \rightarrow0}|z|\log(1/|z|)\lambda(z)=\frac 1 2$, \vspace*{1mm} \\
(ii)$ \;\,\displaystyle l'_{m,\,n}=\frac{1}{2}{-\frac{1}{2} \choose n}{-\frac{1}{2} \choose m}$, \vspace*{1mm} \\
(iii)$\,l'_{m,\,n}=l'_{n,\,m}$.}
\end{theorema}

\vspace*{5mm}
\hspace*{-17pt}\textbf{Acknowledgement.} I would like to thank Prof. Toshiyoki Sugawa for his helpful comments, suggestion and encouragement. I also want to thank Rintaro Ohno for his considerable time reading through each draft.\\

\providecommand{\bysame}{\leavevmode\hbox to3em{\hrulefill}\thinspace}
\providecommand{\MR}{\relax\ifhmode\unskip\space\fi MR }
\providecommand{\MRhref}[2]{%
  \href{http://www.ams.org/mathscinet-getitem?mr=#1}{#2}
}
\providecommand{\href}[2]{#2}

\newpage
\end{document}